\documentclass[12pt]{article}

\usepackage{amsmath,amssymb,amsfonts,theorem,makeidx,latexsym,epsfig,,subfigure}

\newtheorem{defn}{Definition}[section]

\newtheorem{lemma}[defn]{Lemma}

{\theorembodyfont{\rmfamily}

	\newtheorem{ex}[defn]{Example}}

\newtheorem{thm}[defn]{Theorem}

\newtheorem{prop}[defn]{Proposition}

\newtheorem{cor}[defn]{Corollary}

\numberwithin{equation}{section}

\newcommand{\h}{{\cal H}}

\newcommand{\lk}{\lambda_k}

\newcommand{\ltr}{ L^2(\mathbb R) }

\newcommand{\ltn}{{\ell}^2(\mathbb N)}

\newcommand{\cn}{\mc^n}

\newcommand{\si}{S^{-1}}

\newcommand{\afh}{ \forall f \in \h}

\newcommand{\mn}{\mathbb N}

\newcommand{\mr}{\mathbb R}

\newcommand{\mz}{\mathbb Z}

\newcommand{\mq}{\mathbb Q}

\newcommand{\mts}{ \{E_{mb}T_{na}g \}_{m,n \in \mz}}

\def\bp{{\noindent\bf Proof. \ }}

\def\ep{\hfill$\square$\par\bigskip}

\def\bqs{\begin{equation}}

	\def\eqs{\tag*{$\square$}\end{equation}\par\bigskip}

\def\la{\langle}

\def\ra{\rangle}

\def\ga{\gamma}

\def\ftk{\{f_k\}_{k=1}^\infty}

\def\ctk{\{c_k\}_{k=1}^\infty}

\def\gtk{\{g_k \}_{k=1}^\infty}

\def\etk{\{e_k\}_{k=1}^\infty}

\def\suk{\sum_{k=1}^\infty}

\def\nl{\left|\left|}

\def\nr{\right|\right|}

\def\span{\overline{\text{span}}}

\def\Span{\text{span}}

\def\vn{\vspace{.1in}\noindent}

\def\bop{\begin{op}\rm}

	\def\eop{\end{op}}

\def\cn{{\cal N}}

\def\bee{\begin{eqnarray}}

	\def\ene{\end{eqnarray}}

\def\bes{\begin{eqnarray*}}

	\def\ens{\end{eqnarray*}}

\def\bei{\begin{itemize}}

	\def\eni{\end{itemize}}

\def\bt{\begin{thm}}

	\def\et{\end{thm}}

\def\bc{\begin{cor}}

	\def\ec{\end{cor}}

\def\bpr{\begin{prop}}

	\def\epr{\end{prop}}

\def\bl{\begin{lemma}}

	\def\el{\end{lemma}}

\def\bd{\begin{defn}}

	\def\ed{\end{defn}}

\def\bex{\begin{ex}}

	\def\enx{\end{ex}}

\def\bfi{\begin{fig}}

	\def\efi{\end{fig}}

\def\inr{\int_{-\infty}^\infty}

\def\lk{\lambda_k}

\def\mno{\mn_0}

\def\mno{\mn_0}
\def\orbit{\{T^n \varphi\}_{n=0}^\infty}
\def\orf{\{T^n f_1\}_{n=0}^\infty}
\def\fkz{\{f_k\}_{k\in \mz}}
\def\biorbit{\{T^k f_0\}_{k\in \mz}}

\def\lkn{\{\lambda_k\}_{k=1}^\infty}
\def\fttk{\{\widetilde{f_k}\}_{k=1}^\infty}

\title{A survey on frame representations via dynamical sampling}

\date{\today}

\author{Ole Christensen, Marzieh Hasannasab }

\begin{document}

\author{Ole Christensen and Marzieh Hasannasab}
%
%
\maketitle

\begin{abstract}
	Dynamical sampling deals  with representations of a frame $\ftk$ as an orbit  $\orbit$ of a linear
and possibly bounded operator $T$ acting on the underlying Hilbert space.  It is known that the desire of boundedness
of the operator $T$  puts severe restrictions on the frame $\ftk.$  The purpose of the paper is to present an overview
of the results in the literature and also discuss various alternative ways of representing a frame; in particular  the class of considered frames can be enlarged drastically by allowing representations using only a subset 
$\{ T^{\alpha(k)} \varphi\}^\infty_{k=1}$ of 
the operator orbit $\orbit$. In general it is  difficult  to specify appropriate values for the scalars $\alpha(k)$ and  
the vector $\varphi;$ however, by accepting an arbitrarily small  and controllable deviation between the given frame $\ftk$ and 
$\{ T^{\alpha(k)} \varphi\}_{k=1}^\infty$  we will be able to do so.
\end{abstract}

\section{Introduction}

The purpose of this paper is to give an overview of the research topic {\it dynamical sampling,}
seen from the frame perspective.  We will also connect dynamical sampling with certain developments
within linear dynamics, e.g., hypercyclic operators and operator orbits.

Dynamical sampling deals with representations of a frame $\ftk$ as an orbit $\orbit$  of a
bounded  linear operator
$T$ acting on the underlying Hilbert space.  We will highlight a number of necessary conditions on
the frame $\ftk$ for such a representation to exist, restrictions that unfortunately exclude most
classical frames.  Motivated by this we discuss a number of
alternative operator representations of frames, e.g., using unions of orbits. Most importantly, we will
consider frame representations using 
suborbits  $\{T^{\alpha(k)}\varphi\}_{k=1}^\infty$  for certain integers $\alpha(k), k\in \mn,$ a step that
enlarge the class of relevant frames considerably.   As a final step we will consider approximate representations
using suborbits 
and show that by allowing a controllable and arbitrarily small deviation between $\ftk$ and $\{T^{\alpha(k)}\varphi\}_{k=1}^\infty$ 
we will be able to specify appropriate choices of  the scalars
$\alpha(k), \, k\in \mn.$

The paper is presented in a problem-driven way, e.g., using Gabor frames and shift-invariant frames
as motivation for the development.  In the rest of the introduction we set the stage by providing the necessary
background from operator theory and frame theory. The material about (exact) operator representations is collected in Section \ref{210623a}, while various methods for approximate representations are considered
in Section \ref{210623b}. In the entire paper $\h$ denotes a separable infinite-dimensional Hilbert space.

\subsection{Operator representation of sequences}

While our main interest is to analyse frames, several of our operator theoretical tools also occur
in a more general context in the literature. For this reason we formulate the following
definitions for general sequences in Hilbert spaces.

\bd Consider an ordered sequence $\ftk= \{f_1, f_2, f_3, \dots\}$ of elements in $\h.$
If there exists a linear operator $T: \Span \ftk \to \h$ such that $f_{k+1}=Tf_k$ for
all $k\in \mn,$ we say that the sequence $\ftk$ is represented by the operator $T.$ \ed

Note that the sequence $\ftk$ being represented by $T$ precisely means that
\bes \ftk=  \{f_1, f_2, f_3,  \dots\}= \{f_1, Tf_1, T^2f_1, \dots\}= \orf.\ens

In operator theoretical terms, this is phrased by saying that the sequence
$\ftk$ is an {\it orbit} of the operator $T$:

\bd  Given a linear operator $T:  \h \to \h$ and some
$\varphi \in \h,$ the sequence $\orbit$
is called the orbit of $\varphi$. \ed

In general,  the ordering of a given sequence $\ftk$ is crucial for
the question whether or not it can be represented by a linear operator.  For example, if $f$ is an arbitrary nonzero vector in $\h,$ the sequence
$\{f,0,0,0, \dots \}$ is represented by the zero-operator, while the sequence
$\{0,f,0,0,\dots\}$ can not be represented by a linear operator.
We also observe that any linearly independent family $\ftk$ can be represented
by a linear operator $T: \Span \ftk \to \h ,$ regardless of how the sequence is ordered; indeed, we simply define the operator $T$ by
$Tf_k:=f_{k+1}, \, k\in \mn,$ and extend it by linearity.  It was shown in
\cite{olemmaarzieh} that for sequences $\ftk$ spanning an infinite-dimensional space,
linear independence is actually equivalent with the property of being representable by
a linear operator:

\bpr \label{511p} Consider any sequence $\ftk$ in $\h$ for which
$\mbox{span}\ftk$ is infinite-dimensional. Then the following are equivalent:
\bei
\item[(i)] $\ftk$ is linearly independent.
\item[(ii)] There exists a linear operator $T: \mbox{span} \ftk \to \h$
such that $\ftk= \orf.$\eni \epr

In the current paper we will exclusively be interested in representability via
{\it bounded} operators. Even for linearly independent sequences, the possibility of
such a representation depends on the ordering, as the next example demonstrates.

\bex  \label{210806b} Let $\etk$ denote an orthonormal basis for $\h,$ and consider the sequence
\bes \ftk = \{e_1, 2e_2, 3e_3, \dots \}= \{ke_k\}_{k=1}^\infty.\ens
Then $\ftk$ is represented by the linear operator $T$ defined by  \bes Te_k= \frac{k+1}{k}\,e_{k+1}, \,
k\in \mn;\ens note that $T$ extends to a bounded operator on $\h,$ with $||T||=2.$
On the other hand, the following reordering of the elements
\bes \gtk= \{e_1, 10e_{10}, 2e_2, 10^2e_{10^2}, 3e_3, 10^3e_{10^3}, \dots, 9 e_9, 10^9 e_{10^9}, 11 e_{11}, 10^{10}e_{10^{10}}, \dots\} \ens
clearly can not be represented by a bounded operator.
 \enx

In the rest of the paper the focus will be on operator representations of frames. Readers who are interested in
corresponding results for general sequences (in Banach spaces) are referred to \cite{olemmaarziehsteidl}.

\subsection{Hypercyclic operators} \label{212907a}
The main focus of the paper is to consider exact and approximate representations of  frames.  Parts of
the flow of the paper is motivated by the fact that there exist very special operators which have
orbits that are    dense  in the underlying Hilbert space and thus  allow to approximate {\it arbitrary} sequences.

\bd A linear operator $T: \h \to \h$ is hypercyclic if there exists $\varphi \in \h$
such that the orbit $\orbit$ is dense in $\h.$ In
this case the vector $\varphi$ is called a hypercyclic vector. \ed

Note that if $T: \h \to \h$ is a hypercyclic operator with hypercyclic vector
$\varphi,$ then the set $\{T^n \varphi\}_{n=N}^\infty=
\{T^n\left(T^N \varphi\right)\}_{n=0}^\infty$ is dense in $\h$ for
all $N\in \mn;$  this implies that $T^N\varphi$ is a hypercyclic vector for
all $N\in \mn,$ and thus that the set of hypercyclic vectors associated with $T$
is dense in $\h.$

The first example of a hypercyclic operator on a separable Hilbert space was
constructed by Rolewicz  \cite{rolewicz}  in 1969. It deals with the {\it  left-shift operator}
$L(x_1,x_2,x_3,\dots):= (x_2,x_3, \dots)$ on $\ltn$:
\bl \label{70911b} Assume that $a > 1 $ and let $L$ denote the left-shift operator on $\ltn.$ Then $T:=aL$ is a hypercyclic operator on $\ltn.$\el

More  recently  hypercyclic operators have also been constructed on $\ltr$:

\bex \label{231018a} In \cite{Chen} it was proved that for any $a\in \mr$ there exists  functions $\omega: \mr \to \mr$  such that the operator $Tf(x):=\omega(x)f(x-a)$
is hypercyclic on $\ltr$. In particular, for the case  $a=0$  we can take any  continuous decreasing function
such that
\bes \omega(x)= \begin{cases} 2 \, &\mbox{if} \, \, x\le 0 \\
1/2  \, \, &\mbox{if} \, x\ge 1.  \end{cases} \ens
\enx

 For a survey of  other  hypercyclic operators we refer to \cite{erdmannsurvey}.
While the existence of hypercyclic operators is surprising in itself, the theory
also contains several beautiful and intriguing results. We collect a few of
them below, and refer to the literature \cite{BM,erdmann} for more results.

\bt \label{210206a} Consider a bounded linear operator $T: \h \to \h.$  Then the following hold:

\bei
\item[(i)]  {\bf  \cite{ansari}}  If $T$ is hypercyclic, then $T^N$ is hypercyclic
for every $N\in \mn;$ moreover, $T$ and $T^N$ have the same set of hypercyclic vectors.
\item[(ii)]  \ \vspace{.05in}   \ {\bf  \cite{Costakis,Peris}}  If there exist a finite collection
of vectors $\varphi_1, \dots, \varphi_n\in \h$ such that the union of their orbits
is dense in $\h,$ then one of the orbits is dense in $\h$ by itself, i.e., $T$ is
hypercyclic.
\item[(iii)]  \  \vspace{.1in} \  {\bf  \cite{Bourdon}} If an orbit $\orbit$
is somewhere dense in $\h,$ it is automatically dense in $\h,$ i.e., $T$ is hypercyclic
with hypercyclic vector $\varphi.$
\eni \et

Note in particular the result stated in
Theorem \ref{210206a} (i):  phrased in words, it says that if $T$ is hypercyclic with hypercyclic
vector $\varphi,$ then the subset  $\{\varphi, T^N \varphi, T^{2N} \varphi, \dots\}$
is also dense in $\h$, for every $N\in \mn!$

Let us conclude this section with  an interesting connection between hypercyclic operators and the celebrated {\it invariant subspace problem,}  a classical open problem that  can be phrased as follows:

\vn{\bf Q:}  Does there exist a bounded linear operator
$T: \h \to \h$ having no closed nontrivial invariant subspace?

\vn   Note that if $T: \h \to \h$ is a bounded linear operator, the closed (sub)space
$V:= \span\{T^n \varphi\}_{n=0}^\infty$ is invariant for any $\varphi \in \h.$ Thus, a
hypercyclic operator for which every nonzero vector is  hypercyclic
would provide evidence for the existence of the type of operator asked for in the invariant
subspace problem; indeed, such an operator is the {\it only} candidate.  Taking Theorem \ref{210206a} (iii) into account
and   following the
formulation by Bourdon and Feldman in \cite{Bourdon}, we can phrase the invariant
subspace problem as follows:

\vn{\bf Q:}  Does there exists a bounded linear operator $T: \h \to \h$ having a nontrivial nowhere dense orbit?

\subsection{Operator representations via suborbits}
We already considered subsets of operator orbits in Section \ref{212907a}.
 This motivates the next definition.

\bd   Consider a linear operator $T:  \h \to \h$ and some
$\varphi \in \h.$   Given an increasing sequence of nonnegative integers,
\bes 0 \le \alpha(1)  < \alpha(2) <  \alpha (3) < \dots, \ens  the sequence $\{T^{\alpha(k)}  \varphi\}_{k=1}^\infty$
is called  a  suborbit of $\orbit$. \ed

The following deep theorem by Halperin et al.  \cite{kitai} shows that we gain a significant freedom
by asking for sequences to be represented by suborbits of bounded operators rather
than the full orbit.

\bt  \label{210623c} Consider
	any linearly independent sequence $\ftk$ in a Hilbert space $\h.$ Then there exists a bounded operator
	$T: \h \to \h$ and an increasing sequence of nonnegative integers
	\bes 0 \le \alpha(1) < \alpha(2) < \alpha(3) <  \cdots\ens such that
	\bee \label{80513a} \ftk= \{T^{\alpha(k)}f_1\}_{k=1}^\infty.\ene  \et

Theorem \ref{210623c} is very appealing as a theoretical result, but it does not explain how one can
choose the operator $T$ and the integers $\alpha(k)$ in practice.  The proof in \cite{kitai} involves the
selection of a dense linearly independent sequence of elements in the underlying Hilbert space $\h$
and does not give a feasible way of identifying $T$ and $\alpha(k).$ For example, if $\etk$ is an ortnonormal basis for $\h,$ it is easy to
see directly that $\etk= \{T^n e_1\}_{n=0}^\infty,$ where the operator $T$ is defined by $Te_k:=e_{k+1},
\, k\in \mn;$ however, even in this simple case the proof of Theorem \ref{210623c}  yields a very complicated
procedure. One of the main goals of the current paper is to provide more applicable and explicit methods for
representing sequences $\ftk$ in an exact or approximate sense, with an explicit given operator $T$ and
direct access to the integers $\alpha(k).$

\subsection{Basic frame theory}  \label{70331p}
In this section we will give a very short   presentation of the key elements in frame
theory.  Only aspects that are directly relevant for the current paper will be addressed; for more information we refer to
\cite{CB} and \cite{Heil}.

A sequence $\ftk$
in  $\h$ is a {\it frame}
for  $\h$ if there exist constants $A,B>0$ such that
\bes A \,||f||^2 \le \suk | \la f, f_k \ra |^2 \le B\, ||f||^2, \, \afh;\ens it is {\it tight} if we
can choose $A=B,$ and
it is a {\it Bessel sequence} if at least the upper frame condition holds. A {\it Riesz basis} is a frame
which is at the same time a basis; alternatively, a Riesz basis is a frame $\ftk$ for which
\bes  \suk c_k f_k=0, \, \ctk \in \ltn \Rightarrow c_k=0, \, \forall k\in \mn.\ens

If $\ftk$ is a Bessel sequence, the {\it synthesis operator} is defined by
\bee \label{60811a} U: \ltn \to \h, \, U\ctk := \suk c_k f_k;\ene
it is well known that $U$ is well-defined and bounded.
An important role is played by the kernel
of the operator $U,$ i.e., the subset of $\ltn$ given by
\bee \label{60811f} \cn(U)=\left\{\ctk \in \ltn ~\bigg|~\suk c_kf_k=0\right\}.\ene

The  {\it frame operator} is defined by
\bes S: \h \to \h, \,  Sf:=UU^*f= \suk \la f,f_k\ra f_k.\ens
For a frame $\ftk,$ the frame operator is bounded, bijective, and self-adjoint; these
properties immediately lead to the important {\it frame decomposition}
\bee \label{70819d} f= SS^{-1}f= \suk \la f, \si f_k\ra f_k, \quad \forall f\in \h.\ene
The sequence $\{\si f_k\}_{k=1}^\infty$ is also a frame; it is called the {\it canonical dual frame.}
Frames that are {\it not} Riesz bases are   {\it overcomplete,} and
there exists $\gtk \neq \{\si f_k\}_{k=1}^\infty$
such that
\bee \label{70819f}  f= \suk \la f,g_k\ra f_k, \quad \forall f\in \h.\ene

Any frame $\gtk$ satisfying \eqref{70819f} for a given frame $\ftk$ is called a
{\it dual frame} of $\ftk.$

The {\it excess} of a frame $\ftk$ is the maximal number of elements that can be removed yet leaving a frame.
It is well-known that the excess equals $\mbox{dim}\, \cn(U),$  the dimension of the kernel of the
synthesis operator; see  \cite{Balan}.

\subsection{Structured function systems } \label{210611c}
In this section we will give a short introduction to a number of explicit frames in the Hilbert space $\ltr.$ These
frames will illustrate and motivate several results throughout the paper. As for Section \ref{70331p} only results with direct relevance to the current paper are discussed.

In order to introduce the relevant frames  we first need to consider a number of operators on $\ltr.$
For $a\in \mr,$ define the {\it translation operator}
\bee \label{210611a} T_a: \ltr \to \ltr, \, T_af(x):= f(x-a)\ene and the {\it modulation operator}
\bee \label{210611b} E_a: \ltr \to \ltr, \,  E_af(x):=e^{2\pi i ax} f(x).\ene  The
translation operators and the modulation operators are unitary.

Given a function $\varphi\in \ltr$
and some $b>0$,  the associated {\it shift-invariant system} is given by
$\{T_{kb}\varphi\}_{k\in\mz}.$ The frame properties of such systems are well understood, see, e.g.,  \cite{BL1,CDH,Balan}.
In particular, regardless of the choice of $\varphi\in \ltr$ and the parameter  $b>0,$ the system can not be a frame
for all of $\ltr$ but only for the subspace $\span \{T_{kb}\varphi\}_{k\in\mz}.$ It is also known
that $\{T_{kb}\varphi\}_{k\in\mz}$ is linearly independent for all $\varphi \neq 0.$

Classical examples of  frames of the form $\{T_{kb}\varphi\}_{k\in\mz}$ are obtained by taking $\varphi$
to be the sinc-function,
\bes \mbox{sinc}(x):= \begin{cases}  \frac{\sin x}{x} &\mbox{if} \, x\neq 0, \\
1   &\mbox{if} \, x= 0. \end{cases}  \ens
Let us explain this in more detail. First, define  the Fourier transform of $f\in L^1(\mr)$ by \bes \widehat{f}(\ga):=
\inr f(x)e^{-2\pi i \ga x}dx;\ens we extend the Fourier transform in the standard way to a unitary operator
on $\ltr.$ Then  the functions  $\{ T_k  \mbox{sinc} \}_{k\in \mz}$  form an orthnormal basis for
the {\it Paley-Wiener space,}
 \bes PW:=\left\{ f\in \ltr \, \big| \, \mbox{supp} \widehat{f} \subseteq [- \frac12, \frac12 ]\right\}.\ens  Also, the system $\{ T_{k/2}  \mbox{sinc} \}_{k\in \mz}$ can be considered as
the union of two orthonormal bases for $PW,$ and hence form a tight overcomplete frame for $PW.$

It is known that  a finite union of shift-invariant systems at most can form a frame for a subspace for $\ltr.$ 
However, by acting on a shift-invariant system with an infinite number of
modulation operators, it is possible to construct frames for $\ltr.$ This leads to the so-called
{\it Gabor frames.} More formaully, given  some $a,b>0$ and a function
$g\in \ltr$, the associated  {\it Gabor system} is the collection of functions
given by
\bes \mts = \{e^{2\pi imbx}g(x-na)\}_{m,n\in \mz}.\ens
The following result collects the information about Gabor systems that is needed in the current paper;
much more information can be found, e.g., in \cite{G2,CB}.

\bpr \label{61213a} Let $g\in \ltr \setminus \{0\}.$ Then the following hold:

\bei
\item[(i)] $\mts$ is linearly independent.
\item[(ii)]  \ \vspace{.05in} \ If $\mts$ is a frame for $\ltr,$ then $ab\le 1.$
\item[(iii)] \ \vspace{.05in} \ If $\mts$ is a frame for $\ltr,$ then $\mts$ is a Riesz basis if and only if $ab= 1.$
\item[(iv)] \ \vspace{.05in} \  If $\mts$ is an overcomplete frame for $\ltr$, then it has infinite excess.
\eni
\epr
The result in (i) was proved in \cite{Linn}  (hereby confirming a conjecture stated in
\cite{HRT});  (ii) \& (iii) are classical results \cite{G2,CB}, and (iv) is proved in \cite{Balan}.

\section{Frames and dynamical sampling} \label{210623a}

\subsection{Frame properties of operator orbits} \label{210624b}

The analysis of the frame aspects related to dynamical sampling was initiated in the papers
\cite{A1,A2}  by Aldroubi et al. The paper \cite{A2} gives a complete treatment of the finite-dimensional
case; we will not go into this topic as we only deal with the infinite-dimensional case in the current paper.

The approach in \cite{A1,A2} puts operator theory at the central spot.  The central questions can be
formulated as follows:

\bei
\item[(i)]  Given a bounded operator $T: \h \to \h$ belonging to a certain class, is it possible to choose
$\varphi \in \h$ such that the orbit $\orbit$ is a basis?
\item[(ii)]   \ \vspace{.05in} \ Given a bounded operator $T: \h \to \h$ belonging to a certain class, is it possible to choose
$\varphi \in \h$ such that the orbit $\orbit$ is a frame?
\eni
The papers \cite{A1,A2,olemmaarzieh} provide positive answers to these questions in certain cases,
to be described in detail in Sections \ref{210406a}  \& \ref{210406b}. However, except for these cases,
most results in the literature are negative.  We state a number of such results next.

\bpr Consider a bounded opearor $T: \h \to \h.$  Then the following hold:

\bei \item[(i)] {\bf \cite{A2}}  If $T$ is normal,
	then $\{T^n \varphi \}_{n=0}^\infty$ is not a basis.
	\item[(ii)]  \ \vspace{.05in} \ {\bf \cite{A3} } If $T$ is unitary, then $\{T^n \varphi \}_{n=0}^\infty$
	is not a frame.
	\item[(iii)]  \ \vspace{.05in} \ {\bf \cite{olemarzieh-3}} If $T$ is compact, then $\{T^n \varphi \}_{n=0}^\infty$
	is not a frame.

	\eni

\epr

If $T$ is a hypercyclic operator,  clearly $\{ T^n\varphi\}_{n=0}^\infty$ can not be a frame for any
hypercyclic vector $\varphi\in \h.$   We will now prove that $\{ T^n  \eta\}_{n=0}^\infty$ can not be a frame for {\it  any} choice of
$\eta \in \h.$  For this purpose we need the following result by Aldroubi and Petrosyan \cite{A3}.

\bl \label{210624a}  Assume that $\orbit$ is a frame for a bounded operator $T: \h \to \h$ and
some $\varphi \in \h.$ Then $(T^n)^*f\to 0$ as $n\to \infty,$ for all $f\in \h.$  \el

\bp Fix any $f\in \h,$ and let $k\in \mn.$ Then, denoting a lower frame bound for $\orbit$ by $A,$  we have that
\bes A\, ||(T^k)^*f||^2 \le \sum_{n=0}^\infty | \la (T^k)^*f, T^n \varphi \ra|^2 =\sum_{n=k}^\infty | \la f, T^n \varphi \ra|^2 \to 0 \, \mbox{as} \, k\to \infty.\ens
Thus also $(T^k)^*f\to 0$ as $k\to \infty,$
 \ep

\bpr \label{210806d} Assume that $T: \h\to \h$ is hypercyclic. Then $\{T^n \eta\}_{n=0}^\infty$
and $\{(T^*)^n \eta\}_{n=0}^\infty$ can not be frames for
any choice of $\eta \in \h.$ \epr

\bp Let $\varphi\in \h$ be a hypercyclic vector and consider any $\eta\in \h \setminus \{0\}.$ Then the scalar sequence $\{\la \eta, T^n \varphi\ra\}_{n=0}^\infty =
\{ \la (T^*)^n \eta, \varphi\ra\}_{n=1}^\infty$ is unbounded, which implies that the sequence of norms
$\{||(T^*)^n \eta||\}_{n=1}^\infty$ is unbounded. Using Lemma \ref{210624a} it now follows that the
sequence $\{T^n \eta\}_{n=0}^\infty$ is not a frame.
For the proof of the second claim, considering again any $\eta \neq 0$  and still letting $\varphi$ denote a
hypercyclic vector,
\bes \sum_{n=0}^\infty | \la \varphi, (T^*)^n \eta\ra|^2=
\sum_{n=0}^\infty | \la T^n\varphi, \eta\ra|^2=\infty,\ens
which implies that $\{(T^*)^n \eta\}_{n=0}^\infty$ is not a frame. \ep

Yet we have not seen any examples of a frame with an operator representation  $\orbit.$  However, let us prove already now that if $\orbit$
is a frame for a surjective and bounded operator $T,$  the frame property is preserved if an arbitrary finite
number of elements is removed from the frame:

\bpr \label{210614a} Assume that $\orbit$ is a frame for $\h$ for some bounded surjective operator $T: \h \to \h$ and
some $\varphi\in \h.$   Then $\{T^n \varphi\}_{n \in \mno \setminus I}$ is a frame for $\h$  for an arbitrary
finite index set $I \subset \mno.$   In particular, $\orbit$ has infinite excess.  \epr

\bp Consider any $N\in \mn.$  By removing the first $N$ elements of the frame $\orbit$ we are left with the
family $\{T^n \varphi\}_{n=N}^\infty=  \{T^N T^n \varphi\}_{n=0}^\infty;$  since the operator $T^N$ is
bounded and surjective, this family is a frame.   \ep

\subsection{The Carleson frame} \label{210406a}

The first construction of a frame
of the form $\{T^n \varphi \}_{n=0}^\infty$ was obtained by Aldroubi et all in \cite{A1} and further discussed in \cite{A2,A3}.  We
will formulate the result in our setting of an arbitrary separable infinite-dimensional Hilbert space.

\bt \label{70422a}   Let $\etk$ be an orthonormal basis for $\h$  and consider a bounded operator $T: \h \to \h$ such that
$Te_k= \lambda_k e_k$ for a bounded sequence of complex scalars $\lambda_k.$ Also, let
$\varphi\in \h.$ Then  $\{T^n\varphi \}_{n=0}^\infty$ is a frame for $\h$ if and only if the following conditions
are satisfied:
\bei
\item[(i)] $|\lambda_k| <1 $ for all $k\in \mn;$
\item[(ii)]  \ \vspace{.05in} \ $|\lambda_k| \to 1$ as $k\to \infty;$
\item[(iii)]  \ \vspace{.05in} \ The sequence $\{\lambda_k\}_{k=1}^\infty$ satisfies the Carleson condition,
i.e.,
\bee\label{carleson}\displaystyle\inf_n \prod_{n\neq k}\frac{|\lk-\lambda_n|}{|1-\lk\overline{\lambda_n}|}>0;\ene
\item[(iv)] \ \vspace{.05in} \ $\varphi=  \sum_{k=1}^\infty m_k\sqrt{1- |\lambda_k|^2}e_k$ for a scalar-sequence
$\{m_k\}_{k=1}^\infty$ that is bounded below away from zero and above.
\eni
\et

For the sake of easy reference we will call the frames arising from the conditions in Theorem \ref{70422a} for
{\it Carleson frames.} The following result 
(see, e.g., \cite[Thm. 9.2]{duren70}) gives an easy verifiable criterion for the Carleson condition
\eqref{carleson} to hold.
\bpr \label{1406a}
	Let $\lkn$ be a sequence of distinct complex numbers such that $\lk| <1 $ for all $k\in \mn.$ If
there exists $c\in ]0,1[$ such that
	\bee \label{1606c}  \frac{1-|\lambda_{k+1}|}{1-|\lambda_{k}|}\leq c<1, \quad \forall k\in\mn,\ene
	then $\lkn$ satisfies the Carleson condition.
	If $\lkn$ is positive and increasing,   the condition \eqref{1606c} is also necessary for $\lkn$ to satisfy the Carleson condition.
\epr

Based on Proposition \ref{1406a} it is easy to give an explict example of a sequence satisfying the Carleson
condition:
\bc For every $\alpha>1$, the sequence $\lkn=\{1-\alpha^{-k}\}_{k=1}^\infty$ satisfies  the Carleson condition.
\ec

Carleson frames  have a number of very special features.  For example, under a very weak condition, arbitrary
finite subsets  can be removed without destroying the frame property:

\bc  Consider a Carleson frame $\orbit$ as in Theorem   \ref{70422a} and assume that $\lambda_k\neq 0$
for all $k\in \mn.$  Then $\{T^n \varphi\}_{n=N}^\infty$ is a frame for $\h$ for any $N\in \mn.$ \ec

\bp Since $\lambda_k\neq 0$ for all $k\in \mn,$ the condition (ii) in Theorem \ref{70422a} implies
that there is a $C>0$ such that $|\lambda_k| \ge C$ for all $k\in \mn;$ hence the operator $T$ in
Theorem \ref{70422a} is surjective.  The result now follows from Proposition \ref{210614a}. \ep

\subsection{Frames and orbit representations}  \label{210406b}

We will now continue the theme from Section \ref{210624b}, but with an important change of focus. While
Section \ref{210624b} was putting the operator $T$ in the central spot and asking for frame properties
of the associated orbits, we will now take a frame $\ftk$ as the starting point and analyze when and how
it can be represented as an orbit of a bounded operator.  Recall from Proposition \ref{511p} that a representation
as the orbit of a {\it linear} operator is possible if $\ftk$ is linearly independent; thus, the key issue is to
determine when the operator $T$ can be chosen to be bounded.

We first state a result appearing in \cite{olemarzieh-3}.
We will need the {\it  right-shift operator} on $\ltn,$ defined by
\bee \label{70606a} {\cal T}: \ltn \to \ltn, {\cal T}\ctk= \{0, c_1, c_2, \dots\}.\ene

	\bt\label{70510a}
Consider a  frame  $\ftk$ with frame bounds $A,B.$ Then the following are equivalent:

\bei
\item[(i)] The frame has a representation $\ftk=\{T^nf_1\}_{n=0}^\infty$ for some
bounded operator $T:\h \to \h. $
\item[(ii)] \ \vspace{.05in} \ For some dual frame $\gtk$ (and hence all),
\bee \label{70510b} f_{j+1}= \suk \la f_j, g_k\ra f_{k+1}, \, \forall j\in \mn.\ene

\item[(iii)] \ \vspace{.05in} \ The kernel $\cn(U)$
  of the synthesis operator $U$ is invariant  under the right-shift operator ${\cal T}.$
  \eni In the affirmative case, letting $\gtk$ denote an arbitrary dual frame of $\ftk,$
the operator $T$ has the form
\bee \label{70505b}  Tf= \sum_{k=1}^\infty \la f, g_k\ra f_{k+1}, \, \forall f\in \h,\ene
and $1 \le \| T\|\le \sqrt{BA^{-1}}.$\et

Note that if $\ftk$ is a Riesz basis, then $\ftk$ and the (unique) dual frame $\gtk$ form a biorthogonal system, and
hence the condition \eqref{70510b} is trivially satisfied. Thus Theorem \ref{70510a} immediately shows
that every Riesz basis is the orbit of a bounded operator, a result that  was first proved directly in \cite{olemmaarzieh}:

\bc \label{210624c}   If $\ftk$ is a Riesz basis, there exists a bounded operator $T: \h \to \h$ such that
$\ftk= \orf.$ \ec

Note that any reordering of a Riesz basis is a Riesz basis itself; thus, no matter how the elements in a Riesz
basis are ordered, we can apply Corollary \ref{210624c}. On the other hand, we saw
in Example \ref{210806b}  that for general linearly independent sequences in Hilbert spaces, the existence of
a bounded operator representation depends on the ordering of the sequence.   For general overcomplete frames,
the question of representability via the orbit of a bounded operator is also extremely sensitive to reorderings of
the elements, as demonstrated by results in \cite{olemarzieh-3}.

The simplicity of the proof of Corollary \ref{210624c}  could make us optimistic with regard to
construction of other classes of frames having an orbit representation via a bounded operator, but
unfortunately this is misleading.  Indeed, the Carleson frames and the Riesz bases are the only explicitly
available examples in the literature of frames having such a representation. In particular,  it was  shown in  \cite{olemarzieh-3} that a so-called {\it near-Riesz basis}  (i.e., a  frame with finite excess)  never can be represented via
a bounded operator:

\bpr Let $\ftk$ denote a frame which is a union of a Riesz basis and a finite non-empty collection
of vectors. Then $\ftk$ can not be represented as the orbit of a bounded operator.\epr
 Another serious restriction  was proved in \cite{CHP}. It excludes all overcomplete shift-invariant frames
and Gabor frames from having a representation as the orbit of a bounded operator:

\bt \label{71218d} Assume that $\orbit$ is an overcomplete frame for
some $\varphi\in \h$ and a bounded operator $T:\h \to \h.$  Then for each $f\in\h$ we have that
\begin{eqnarray} \label{71219d}
T^nf\to 0\quad\text{as }n\to\infty.
\end{eqnarray}
\et

The result in Theorem  \ref{71218d} is very disappointing. 
We have already seen that shift-invariant frames and Gabor frames automatically are linearly independent,
so by Proposition \ref{511p} such systems can be represented as an orbit of a {\it  linear} operator $T;$ however,
in the overcomplete case $T$ is forced to be unbounded. Since the problem is the boundedness and not
the {\it existence } of a operator representation, it is natural to ask whether boundedness can be achieved by
allowing a ``scaling factor" in the operator representation. We formulate it as an open question:

\vn{\bf  Q:} Do there exist overcomplete Gabor frames $\mts$ such
that  an appropriate
ordering $\ftk$ of the frame elements has a representation
\bee \label{70821b} \ftk = \{a_nT^n \varphi\}_{n=0}^\infty, \ene
for some  scalars $a_n \neq0$, a bounded operator $T: \ltr \to \ltr,$ and some $\varphi \in
\ltr$?

\vspace{.1in}
Let us end the section with result connecting frame orbits and shift-invariant subspaces of the Hardy space. Recall that letting
 $\mathbb{D}$ denote the open unit disc in the complex plane, the Hardy space $H^2(\mathbb{D})$ is defined as
\[ H^2(\mathbb{D})=\{ f:\mathbb{D}\to\mathbb{C}:\, f(z)=\sum_{n=0}^{\infty}a_nz^n,\{a_n \}_{n=0}^\infty\in\ell^2(\mn_0)\}.  \]  The Hardy space is a Hilbert space, with the inner product defined by
$$\langle f,g \rangle=\sum_{n=0}^\infty a_n\overline{b}_n, \, \, \mbox{for}  \, \,   f=\sum_{n=0}^\infty a_nz^n, g=\sum_{n=0}^\infty b_nz^n.$$
Furthermore the sequence of functions $\{z^n\}_{n=0}^\infty$ is an orthonormal basis for $H^2(\mathbb{D})$.
The {\it shift operator } on the Hardy space is defined by
\bes {\cal S}:H^2(\mathbb{D})\to H^2(\mathbb{D}),{\cal S}f(z)=zf(z),  \, z\in \mathbb{D}  \ens
Finally let $1_{\mathbb{D}} $ denote the constant function, i.e., $1_{\mathbb{D}}(z) = 1$.

\bt \label{210714a} Consider a bounded operator $T: \h \to \h$ and let $\varphi \in \h.$ Then
 $\orbit$ is a frame  if and only if there exists a shift invariant subspace ${\cal F}\subset H^2(\mathbb D)$  and a bounded bijective operator $W: \cal F^\perp \to\h$ such that
\[
T = W P_{{\cal F}^\perp}{\cal S }W^{-1}, \quad \varphi = W P_{{\cal F}^\perp}1_{\mathbb D}.
\]
In the affirmative case, $\orbit$ is a Riesz basis if and only if ${\cal F}=\{0\}$.
\et

\bp In the entire proof we let $\la \cdot, \cdot \ra$  denote the inner product on $H^2(\mathbb{D}).$  First assume that
 $\orbit$ is a frame.   Consider the linear operator $V:H^2(\mathbb D)\to \h$  defined by
\bee\label{def:operatorV}
Vf = \sum_{n=0}^\infty  \la f, z^n\ra T^n \varphi , \  f\in H^2(\mathbb{D}).
\ene The operator $V$ is bounded; indeed, letting
$B$ denote a Bessel bound for $\orbit$, we have
 \bes
 \| Vf \|^2 = \|  \sum_{n=0}^\infty   \la f, z^n\ra T^n \varphi \|^2 \leq B   \sum_{n=0}^\infty  |  \la f, z^n\ra |^2 = B \| f \|^2.
\ens
The frame assumption also implies that the operator $V$ is surjective. Moreover, for $f \in H^2(\mathbb{D})$ we have
\bee\label{intertwining}
TVf &=& T\sum_{n=0}^\infty   \la f, z^n\ra T^n \varphi = \sum_{n=0}^\infty   \la f, z^n\ra T^{n+1}\varphi\nonumber\\
&= & V \sum_{n=0}^\infty   \la f, z^n\ra z^{n+1} = V{\cal S}f. \label{eq:interchangable}
\ene
This shows in particular that $\ker V$ is a closed shift-invariant subspace of $H^2(\mathbb D)$.
Now let ${\cal F}:=\ker V$ and consider the restriction of the
operator $V$ to the subspace   ${\cal F}^\bot, $ i.e.,  let  $W:= V_{|_{{\cal F}^\perp}}:{\cal F}^\perp \to \h$. Also let $P_{{\cal F}^\perp}: H^2(\mathbb D)\to {\cal F}^\perp $ be the orthogonal projection. Then, using
 \eqref{intertwining}, on the subspace ${\cal F}^\perp,$ we have
\bes  TW=TV= V{\cal S} =   V  P_{{\cal F}^\bot}     {\cal S} =    W  P_{{\cal F}^\bot}     {\cal S}  \ens
Therefore the invertibility of  $W$  yields that
 $T = W P_{{\cal F}^\perp}{\cal S}W^{-1};$
 also $\varphi = V 1_{\mathbb D} =V P_{{\cal F}^\perp} 1_{\mathbb D}=W P_{{\cal F}^\perp} 1_{\mathbb D}$.

Conversely,  let ${\cal F}\subset H^2(\mathbb D)$ be a closed shift-invariant subspace  and $W: {\cal F}^\perp \to \h$ be an invertible bounded operator.  Since $\{{\cal S}^n 1_{\mathbb D}\}_{n=0}^\infty $ is an orthonormal basis for $H^2(\mathbb D)$, the sequence  $\{P_{{\cal F}^\perp}{\cal S}^n 1_{\mathbb D}\}_{n=0}^\infty $ is a frame for ${\cal F}^{\perp}$. Now let $T:= W P_{{\cal F}^\perp} {\cal S}W^{-1}$ and $\varphi := WP_{{\cal F}^\perp} 1_{\mathbb D}$.
 Since ${\cal F}$ is invariant under the action of ${\cal S}$, we have
  \[
  P_{{\cal F}^\perp}{\cal S}  = P_{{\cal F}^\perp}{\cal S} P_{{\cal F}^\perp}
  \]
   which implies that 	$(P_{{\cal F}^\perp}{\cal S})^n  = P_{{\cal F}^\perp}{\cal S}^n$,  for every $n\in \mathbb N_0$. Thus
 \bes
 T^n \varphi &=& W (P_{{\cal F}^\perp} {\cal S})^nW^{-1} \varphi =W (P_{{\cal F}^\perp}{\cal S})^n P_{{\cal F}^\perp}1_{\mathbb D}
  = WP_{{\cal F}^\perp}{\cal S}^n 1_{\mathbb D}
 \ens
 Therefore the collection $\{ T^n\varphi\}_{n=0}^\infty $= $\{WP_{{\cal F}^\perp}{\cal S}^n 1_{\mathbb D}\}_{n=0}^\infty $ is a frame for $\h$.  \ep

Theorem \ref{210714a} is theoretically appealing, but it is not instrumental in terms of constructing more
examples of explicitly given frames $\orbit.$ The difficulty in constructing such frames is the motivation
for the remaining sections in the paper.

\subsection{Frames and bi-orbit representations} \label{210806a}
In this section we will consider a frame, indexed as $\fkz;$ clearly,  we can bring any frame into this form
by a reordering and reindexing.

A seemingly innocent modification of the topic from Section \ref{210406b}  would be
to ask for representations of a given frame $\fkz$ of the form $\biorbit$ for an
{\it invertible}  operator $T: \Span \fkz \to \Span \fkz; $ such a representation is called
a {\it bi-orbit representation.}  Interestingly, this leads to a theory that in certain aspects is very similar
to the operator representations discussed in the previous sections, but also
is completely different at some points. In this section we will highlight
some of the similarities and differences for the two ways of representing a frame.

First, it was proved in \cite{olemmaarzieh-2} that Proposition \ref{511p} has
a completely similar version for bi-orbit representations: any linearly independent
frame $\fkz$ has a representation $\fkz = \{T^k f_0\}_{k\in \mz}$ for some
linear and invertible operator $T: \Span \fkz \to \Span \fkz.$
It was also proved that Theorem \ref{70510a} has a parallel version: indeed, if
additionally  the kernel
of the synthesis operator is invariant under right   {\it and}  left-shifts,  the
operator $T$ extends to a bijective and bounded operator $T: \h \to \h.$

However, when it comes to concrete manifestations of frames of the form $\biorbit$
as well as other properties for such frames, there are remarkable difference between
orbit/bi-orbit representations. This becomes apparent already from the following example.

\bex As discussed in Section \ref{210611c}, the shift-invariant system $\{ T_k  \mbox{sinc} \}_{k\in \mz}$ is an orthnormal basis for
the Paley-Wiener space $PW$ of functions in $\ltr$,  while $\{ T_{k/2}  \mbox{sinc} \}_{k\in \mz}$ is an
overcomplete frame for the same space.  Since  $\{T_{kb}\varphi\}_{k\in \mz}=
\{(T_{b})^k\varphi\}_{k\in \mz},$ such frames are born to have a bi-orbit representation. However,
since $||T_{kb}\varphi||=||\varphi||$ for all $k\in \mz$ it is clear from
Theorem \ref{71218d} that an overcomplete shift-invariant frame can not have a representation
$\orbit$ for a bounded operator $T.$
 \enx

Recall that Lemma \ref{210624a}  and Theorem \ref{71218d} provided necessary conditions for 
an orbit $\orbit$ to be a frame for a bounded operator $T.$  The next result, proved in \cite{CHP}, shows
that the corresponding conditions on a bi-orbit $\{T^n \varphi\}_{n\in \mz}$ are fundamentally different:

\bt \label{210629a} Assume that $\{T^n \varphi\}_{n\in \mz}$ is a  frame with bounds $A,B$ for
some $\varphi\in \h$ and a bounded invertible operator $T:\h \to \h.$  Then for each $f\in\h$ we have that
$$
\|T^nf\|\,\ge\,\sqrt{\frac A B}\,\|f\|\qquad\text{and}\qquad\|(T^*)^{n} f\|\,\ge\,\sqrt{\frac A B}\,\|f\|.
$$
\et

Theorem \ref{210629a} makes it natural to ask whether overcomplete Gabor frames
have representations as bi-orbits of a bounded operator. We phrase the question as stated originally in
\cite{beauty}:

\vspace{.2in}\noindent {\bf Q:} Do there exist  overcomplete Gabor frames $\mts$ such
that  an appropriate
ordering $\fkz$ of the frame elements has a representation
\bes \fkz = \{T^n \varphi\}_{n=-\infty}^\infty, \ens
for some  bounded operator $T: \ltr \to \ltr$ and some $\varphi \in \ltr$?

\vn In \cite{Corso} Corso  provided a partial solution to the question, by showing that the answer is negative
if $ab \notin \mq,$ as well as under certain support conditions on the function $g.$  However, the
general question is still open.

\subsection{Frame representations via multi-orbits}

The difficulty in constructing explicitly given frames $\orbit,$ discussed in detail in Section \ref{210406b}, makes it
natural to explore various modifications of the basic idea. In this section we will discuss the additional flexibility that
is obtained by allowing a frame representation using a {\it finite} number of operator orbits. Let us first connect
the topics in Section \ref{210406b} and Section \ref{210806a} by showing that any frame with a bi-orbit representation
$\{T^n \varphi\}_{n\in \mz}$ can be considered as the union of two orbits:

\bex Consider any frame that has
a representation as a bi-infinite orbit
$\{T^n \varphi\}_{n\in \mz}$ for a bounded
and bijective operator
$T: \h \to \h.$ Then
\bee \label{80923h} \{T^n \varphi\}_{n\in \mz}=
\{T^n \varphi\}_{n=0}^\infty \cup
\{(T^{-1})^n T^{-1}\varphi\}_{n=0}^\infty,\ene
i.e., the frame  $\{T^n \varphi\}_{n\in \mz}$ can be represented as a union of two orbits, generated by the
bounded operators $T$ and $T^{-1}.$
  \enx

The following result, originally proved in \cite{olemmaarzieh}, shows that the idea of multi-orbits significantly
enlarges the class of frames that can be analyzed:  for example, all the classical frames in harmonic analysis
like Gabor frames and wavelet frames have such representations!

\bt \label{50612c} Consider a frame $\ftk$ which is
norm-bounded below. Then there is a finite collection of
vectors from $\ftk,$ to be called $\varphi_1, \dots, \varphi_J,$ and
corresponding bounded operators $T_j: \h \to \h$ with closed range, such that
\bee \label{50614a} \ftk = \bigcup_{j=1}^J \{T_j^n \varphi_j\}_{n=0}^\infty.\ene

\et

While Theorem \ref{50612c} is theoretically appealing, its practical applicability is limited by the fact that 
the proof is based on
the Feichtinger Theorem  and thus does not give direct access to the operators $T_j.$
More results on multi-orbit representations can be found in \cite{CMPP}.

\section{Approximate operator representations} \label{210623b}

In Section \ref{210623a} the focus is on frames $\ftk$  having an exact match with an operator orbit
$\orbit;$ unfortunately,  the results show that this is a very restrictive condition and that the standard
frames in $\ltr$ are excluded.   In the current section we will  consider two relaxations on the orbit 
representations, which in combination yield a theory covering a much larger class of frames.
To be more precise, recall that  Theorem \ref{210623c}  by Halperin et all showed that {\it suborbit representations}
 $\{ T^{\alpha(k)} \varphi \}_{k=1}^\infty$
via bounded operators are available  for any linearly independent family $\ftk;$ however, in general
we do not have a feasible procedure to identify the operator $T$ and the appropriate integers $\alpha(k).$
In this section we will show that for a large class of frames we can obtain {\it approximate suborbit representations,}
where it is possible to specify as well the operator $T$ as the integers $\alpha(k).$  In the entire section we follow
the presentation in \cite{OM-1}; we refer to that paper for proofs and more information.

We begin with a formal definition, explaining the exact meaning of ``approximation." In words, a
frame $\ftk$ is approximated by a sequence $\fttk$ if $\fttk$ satisfies the standard perturbation
condition within frame theory:

\bd \label{191001a} Let $\ftk$ be a frame for $\h.$ Given any $\epsilon >0,$
a sequence $\fttk \subset \h$ is called an $\epsilon$-approximation of
$ \{f_k\}_{k\in I}$ if
\bee \label{151018a}  \nl \sum c_k (f_k - \widetilde{f_k})\nr^2
\le \epsilon \, \sum |c_k|^2\ene
for all finite sequences $\{c_k\}.$\ed

For sufficiently small values of $\epsilon,$ an $\epsilon$-approximation $\fttk$ of a frame $\ftk$
is itself a frame and shares several key properties
of the frame, e.g., its excess. Furthermore, the synthesis
operator and the frame operator for $\ftk$ are
approximated ``well" by the corresponding operators associated with
$\fttk$:

\bt\label{181217ens}  Consider a frame $\ftk$
for $\h$  with frame bounds $A,B$, and assume that
$\fttk \subset \h$ is an $\epsilon$-approximation of
$ \ftk$ for some $\epsilon \in ]0, A[.$
Then the following hold:

\bei\item[(i)]  $ \fttk $ is a frame with bounds $A(1-\sqrt{ \frac{\epsilon}{A}  })^2$ and $B(1+\sqrt{ \frac{\epsilon}{B}  })^2,$ with the same excess as
$ \ftk.$
\item[(ii)] Denoting the synthesis operators and frame operators of $\ftk$ and
    $\fttk$  by $U, \widetilde{U}$, respectively, $S,\widetilde{S}$, we have \[ \| U- \widetilde{U} \| \leq \sqrt{\epsilon},\quad
  \| S - \widetilde{S} \|\leq \sqrt{\epsilon B}  \left(2+\sqrt{ \frac{\epsilon}{B}  }  \right),\]
  and \[
    \| S^{-1} - \widetilde{S}^{-1} \|\leq  \frac{ \sqrt{\epsilon B} (2+\sqrt{ \frac{\epsilon}{B}}) }{A^2(1-\sqrt{ \frac{\epsilon}{A}  })^2}. \]
\eni
\et

The connection to approximate frame representations using suborbits of a bounded operator is
explained in the next result, which provides a natural way to satisfy the condition \eqref{151018a}:

\bc\label{181217en}  Let $\ftk$ be a frame
for $\h$  with lower frame bound $A$. Also,  let  $\varphi \in \h$
and consider a
bounded operator $T: \h \to \h$. Assume
that  for a given $\epsilon \in ]0, A[,$ and for any $k\in \mn$ there exists a nonnegative integer
$\alpha(k)\in\mn_0$ such that
\bee \label{011018p}  \| f_k -  T^{\alpha(k)} \varphi \|^2 \leq\frac{\epsilon}{2^k}.\ene
Then $\{ T^{\alpha(k)} \varphi \}_{k=1}^\infty$ is an $\epsilon$-approximation
of $\ftk$.
\ec

Note that 
Corollary \ref{181217en}  applies to any frame and any  hypercyclic
operator:

\bex \label{171218a} Let $T:\h \to \h$ be a hypercyclic operator with hypercyclic vector $\varphi.$  Then, for any frame
$\ftk$ for $\h$ and any given $\epsilon >0$ there exists
nonnegative integers $\alpha(k), \, k\in \mn,$ such that
\eqref{011018p} holds. Thus, for sufficiently small
values of $\epsilon,$ we obtain a frame $\{T^{\alpha(k)}\varphi\}_{k=1}^\infty$ with the same
excess as the frame $\ftk$ and approximating $\ftk$
in the sense of Definition \ref{191001a}.
 \enx

In general we do not have direct
access to the  numbers  $\alpha(k), \, k\in \mn,$ such that
\eqref{011018p} holds, not even for
hypercyclic operators. However,  under natural conditions we can be very explicit about the choices
of appropriate scalars $\alpha(k)$ for certain choices
of the operator $T.$  In the rest of the section we will use the following

\vn{\bf General setup:}  Fix an orthonormal basis $\etk$ for the Hilbert space $\h,$ and let $\ftk$ denote a frame for $\h$.
 Fix some $\lambda >1$, and define  the {\it scaled
left/right-shift operators} on $\h$  by

\bes \label{231018d}  T \left( \suk c_k e_k \right):= \lambda \suk c_{k+1}e_k, \, \,U \left( \suk c_k e_k \right):= \lambda^{-1} \suk c_{k}e_{k+1}, \, \, \ctk \in \ltn.\ens

For appropriately chosen nonnegative integers $\alpha(k), \, k\in \mn,$
let
\bee \label{011018bn}
\varphi:= \sum_{n=1}^\infty U^{\alpha(n)} f_n.\ene

For the case of finitely supported vectors $\ftk$ (meaning that each vector $f_k$ is a finite linear combination
of elements from $\etk$)
the following result specifies how to choose the powers
$\alpha(k)$ such that the vector $\varphi$ in \eqref{011018bn}
is well-defined and \eqref{011018p} holds:

\bt \label{011018an} Under the conditions in the general setup, assume that $\ftk$  consists of finitely supported vectors
in $\h;$ for
$k\in \mn,$ let $m(k)$ denote the largest index for a nonzero
coordinate in $f_k.$ Let $\{\alpha(k)\}_{k=1}^\infty$ be a
strictly increasing sequence
of nonnegative integers such that $\alpha(1)=0$ and
$\alpha(k+1)-\alpha(k)\ge m(k)$ for all $k\in \mn.$
Then,   for any  $k\in \mn,$
\bee \label{011018cn} || f_k- T^{\alpha(k)}\varphi||^2 \le  \frac{B\lambda^2}{\lambda^2-1} \lambda^{-2[\alpha(k+1)-\alpha(k)]}.\ene
In particular, choosing the nonnegative integers $\alpha(k)$ such that $\alpha(1)=0$ and
for a given $\epsilon \in ]0, A[,$
\bee \label{011018k}
\alpha(k+1)-\alpha(k) \ge \max\left(m(k), \frac{k \,\ln(2)+\ln \left(\frac{B}{\epsilon}\right) + \ln \left( \frac{\lambda^2}{\lambda^2-1}\right)  }{2\ln (\lambda)}\right), \ene
the condition \eqref{011018p} is satisfied, i.e., the conclusions in Theorem \ref{181217ens} hold.

\et

For certain values of the scaling parameter $\lambda$ we can be completely specific about how to choose
the integers $\alpha(k):$

\bc \label{180226a} In the setup of Theorem \ref{011018an}, let $\lambda=\sqrt{2},$ take an upper frame bound of the
form $B=2^N$ for some $N\in \mn$ and a tolerance
$\epsilon =2^{-j}$ for some $j\in \mn.$ Then the following hold.

\bei
\item[(i)] Without any restriction on the support sizes
$m(k)$ of the vectors $f_k,$ the condition \eqref{011018k}
is satisfied if $\alpha(1)=0$ and
\bee \label{061018a} \alpha(k)= (k-1)\left[ N+j + 1+\frac{k}{2}\right] + \sum_{\ell=1}^{k-1} m(\ell), \, \, k\in \mn \setminus\{1\}. \ene\item[(ii)] If
$m(k) \le N+j + 1+k$ for all $k\in \mn,$ the condition \eqref{011018k}
is satisfied if
\bee \label{061018b} \alpha(k)= (k-1)\left[ N+j + 1+\frac{k}{2}\right], \, \, k\in \mn. \ene \eni

\ec

The main condition in Theorem \ref{011018an}, namely, that the frame $\ftk$ is finitely supported,
can be relaxed: indeed, similar but slightly more technical results can be proved for localized
frames. We refer to \cite{OM-1} for the details.  Note also that generalizations to Banach spaces
(for sequences $\ftk$ without assuming any frame property)  have been obtained in \cite{olemmaarziehsteidl}.

%

{\bf \vspace{.1in}
	
	\noindent Ole Christensen\\
	DTU Compute\\
	Technical University of Denmark\\
	Building 303 \\
	2800 Lyngby  \\
	Denmark \\
	Email: ochr@dtu.dk
	
	\vn Marzieh Hasannasab \\
	Institut f\"ur Mathematik\\
	TU Berlin\\
	Stra\ss e des 17. Juni 136\\
	10623 Berlin \\ Germany
	\\
	Email: hasannas@math.tu-berlin.de

}

\end{document}